\documentclass{amsart}
\usepackage{amssymb}

\newtheorem{theorem}{Theorem}
\newtheorem{lemma}{Lemma}
\newtheorem*{corollary}{Corollary}

\newcommand{\C}{\mathbb{C}}

\newcommand{\Z}{\mathbb{Z}}
\newcommand{\N}{\mathbb{N}}
\newcommand{\PP}{\mathbb{P}}
\newcommand{\LL}{\mathcal{L}}

\newcommand{\RR}{\mathcal{R}}
\newcommand{\LLL}{\mathbf{L}}
\newcommand{\OO}{\mathcal{O}}
\newcommand{\ox}{\vec{x}}

\newcommand{\oc}{\vec{c}}
\newcommand{\XX}{\mathbb{X}}

\newcommand{\gf}{\mathfrak{g}}
\newcommand{\kdelta}{\mbox{\boldmath $\delta$}}

\DeclareMathOperator{\rad}{rad}
\DeclareMathOperator{\ind}{ind}
\DeclareMathOperator{\ad}{ad}
\DeclareMathOperator{\GL}{GL}

\DeclareMathOperator{\Coh}{coh}

\DeclareMathOperator{\Hom}{Hom}
\DeclareMathOperator{\Ker}{Ker}

\DeclareMathOperator{\Aut}{Aut}
\DeclareMathOperator{\End}{End}
\DeclareMathOperator{\Ext}{Ext}

\DeclareMathOperator{\rank}{rank}

\DeclareMathOperator{\modu}{mod}
\DeclareMathOperator{\dimv}{\underline{\dim}}

\title{Kac's Theorem for weighted projective lines}
\author{William Crawley-Boevey}
\address{Department of Pure Mathematics, University of Leeds, Leeds LS2 9JT, UK}\email{w.crawley-boevey@leeds.ac.uk}

\begin{document}
\begin{abstract}
We prove an analogue of Kac's Theorem, describing the dimension
types of indecomposable coherent sheaves (or parabolic bundles)
over weighted projective lines in terms of root systems for
loop algebras of Kac-Moody Lie algebras.
We use a theorem of Peng and Xiao to associate a Lie algebra
to the category of coherent sheaves for a weighted projective line
over a finite field, and find elements of this Lie algebra which
satisfy the relations defining the loop algebra.
We use these elements in the proof of our analogue of Kac's Theorem.
\end{abstract}

\subjclass[2000]{Primary 14H60, 16G20.}

\keywords{Weighted projective line, parabolic bundle,
Kac-Moody Lie algebra, loop algebra, Hall algebra.}

\maketitle

\section{Weighted projective lines}
Let $K$ be an algebraically closed field,
let $\PP^1$ be the projective line over $K$,
let $D=(a_1,\dots,a_k)$ be a collection of distinct points of $\PP^1$,
and let $\mathbf{w}=(w_1,\dots,w_k)$ be a \emph{weight sequence}, that is,
a sequence of positive integers.
The triple $\XX = (\PP^1,D,\mathbf{w})$ is called a \emph{weighted projective line}.
Geigle and Lenzing \cite{GL} have associated to each weighted
projective line a category $\Coh\XX$ of coherent sheaves on $\XX$,
which is the quotient category of the category of finitely
generated $\LLL(\mathbf{w})_+$-graded $S(\mathbf{w},D)$-modules,
modulo the Serre subcategory of finite length modules. Here
$\LLL(\mathbf{w})$ is the rank 1 additive group
\[
\LLL(\mathbf{w}) = \langle \ox_1,\dots,\ox_k,\oc \mid w_1 \ox_1 = \dots = w_k \ox_k = \oc \rangle
\]
partially ordered, with positive cone
$\LLL(\mathbf{w})_+ =  \N \oc+\sum_{i=1}^k \N \ox_i$,
and
\[
S(\mathbf{w},D) = K[u,v,x_1,\dots,x_k] / ( x_i^{w_i} - \lambda_i u - \mu_i v ),
\]
with grading $\deg u = \deg v = \oc$ and $\deg x_i = \ox_i$,
where $a_i = [\lambda_i:\mu_i]\in\PP^1$.
Geigle and Lenzing showed that $\Coh\XX$ is a hereditary abelian category
with finite-dimensional Hom and Ext spaces.
The free module gives a structure sheaf $\OO$, and shifting the grading gives
twists $E(\ox)$ for any sheaf $E$ and $\ox\in\LLL(\mathbf{w})$.

Every sheaf is the direct sum
of a `torsion-free' sheaf,
which has a filtration by sheaves of the form $\OO(\ox)$, and a finite-length sheaf,
and the latter are easily described.
There are simple sheaves $S_a$ ($a\in \PP^1\smallsetminus D$)
and $S_{ij}$ ($1\le i\le k$, $0\le j\le w_i-1$).
They have
\[
\dim\Hom(\OO(r\oc),S_{ij})=\kdelta_{j0},
\quad
\dim\Ext^1(S_{ij},\OO(r\oc))=\kdelta_{j1}
\]
where $\kdelta$ is the Kronecker delta function,
and the only extensions between them are
\[
\dim \Ext^1(S_a,S_a) =1,
\quad
\dim \Ext^1(S_{ij},S_{i\ell}) = 1
\quad
(\ell \equiv j-1 \, (\modu \, w_i)).
\]
For each simple sheaf $S$ and $r>0$ there is a unique sheaf $S[r]$, with length $r$
and top $S$, which is \emph{uniserial}, meaning that it has a unique composition series.
These are all the finite-length indecomposable sheaves.

There is a root system associated to $\mathbf{w}$ via the graph
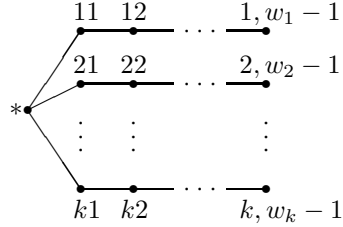
\[
\begin{picture}(110,80)
\put(10,40){\circle*{2.5}}
\put(30,10){\circle*{2.5}}
\put(30,50){\circle*{2.5}}
\put(30,70){\circle*{2.5}}
\put(50,10){\circle*{2.5}}
\put(50,50){\circle*{2.5}}
\put(50,70){\circle*{2.5}}
\put(100,10){\circle*{2.5}}
\put(100,50){\circle*{2.5}}
\put(100,70){\circle*{2.5}}
\put(10,40){\line(2,-3){20}}
\put(10,40){\line(2,1){20}}
\put(10,40){\line(2,3){20}}
\put(30,10){\line(1,0){35}}
\put(30,50){\line(1,0){35}}
\put(30,70){\line(1,0){35}}
\put(85,10){\line(1,0){15}}
\put(85,50){\line(1,0){15}}
\put(85,70){\line(1,0){15}}
\put(70,10){\circle*{1}}
\put(75,10){\circle*{1}}
\put(80,10){\circle*{1}}
\put(70,50){\circle*{1}}
\put(75,50){\circle*{1}}
\put(80,50){\circle*{1}}
\put(70,70){\circle*{1}}
\put(75,70){\circle*{1}}
\put(80,70){\circle*{1}}
\put(30,25){\circle*{1}}
\put(30,30){\circle*{1}}
\put(30,35){\circle*{1}}
\put(50,25){\circle*{1}}
\put(50,30){\circle*{1}}
\put(50,35){\circle*{1}}
\put(100,25){\circle*{1}}
\put(100,30){\circle*{1}}
\put(100,35){\circle*{1}}
\put(3,35.0){*}
\put(27,-1){$k1$}
\put(27,54){$21$}
\put(27,74){$11$}
\put(45,-1){$k2$}
\put(45,54){$22$}
\put(45,74){$12$}
\put(90,-1){$k,w_k-1$}
\put(90,54){$2,w_2-1$}
\put(90,74){$1,w_1-1$}
\end{picture}
\]
whose vertex set $I$ consists of $*$ and
vertices denoted $ij$ or $i,j$ for $1\le i\le k$ and $1\le j\le w_i-1$.
Let $\gf$ be the Kac-Moody Lie algebra (over $\C$)
with generators $e_v,f_v,h_v$ ($v\in I$) and relations
\[
\begin{cases}
[h_u,h_v] =0,
\quad
[e_u,f_v]=\kdelta_{uv} h_v,
\\
[h_u,e_v] = a_{uv} e_v,
\quad
[h_u,f_v] = -a_{uv} f_v,
\\
(\ad e_u)^{1-a_{uv}}(e_v)=0,
\quad
(\ad f_u)^{1-a_{uv}}(f_v)=0
\quad
(\text{if $u\neq v$})
\end{cases}
\]
where the (symmetric) generalized Cartan matrix $(a_{uv})$
has diagonal entries 2 and off-diagonal entries
$-1$ if $u$ and $v$ are joined by an edge and otherwise~0.
The root lattice $\Gamma$ is the free additive group
on symbols $\alpha_v$ ($v\in I$), and there is a symmetric
bilinear form on it defined by $(\alpha_u,\alpha_v) = a_{uv}$.
Now $\gf$ is graded by $\Gamma$,
with $\deg e_v = \alpha_v$, $\deg f_v=-\alpha_v$ and $\deg h_v=0$,
and the root system is
$\Delta = \{ 0\neq \alpha\in \Gamma \mid \gf_\alpha\neq 0 \}$.
Recall that there are real roots,
obtained from the simple roots $\alpha_v$ by a sequence of reflections
$s_u(\alpha) = \alpha-(\alpha,\alpha_u)\alpha_u$,
and there may also be imaginary roots.

The \emph{loop algebra} of $\gf$ is $L\gf = \gf[t,t^{-1}]$,
but more appropriate is an extension $\LL\gf$ with generators
$e_{vr},f_{vr},h_{vr}$ ($v\in I$, $r\in \Z$) and $c$ subject to the relations
\begin{equation}
\label{e:loopalg}
\begin{cases}
\text{$c$ central,}
\quad
[e_{vr},e_{vs}] = 0,
\quad
[f_{vr},f_{vs}] = 0,
\\
[h_{ur},h_{vs}]=ra_{uv}\,\kdelta_{r+s,0}\,c,
\quad
[e_{ur}, f_{vs}] = \kdelta_{uv}\, \left(h_{v,r+s} + r\,\kdelta_{r+s,0}\,c\right),
\\
[h_{ur},e_{vs}] = a_{uv} e_{v,r+s},
\quad
[h_{ur},f_{vs}] = -a_{uv} f_{v,r+s},
\\
(\ad e_{u0})^{1-a_{uv}}(e_{vs}) = 0,
\quad
(\ad f_{u0})^{1-a_{uv}}(f_{vs}) = 0
\quad
(\text{if $u\neq v$}),
\end{cases}
\end{equation}
see~\cite{MRY} and \cite[\S1.3]{Schiffmann}.
The root lattice for either algebra is
$\hat{\Gamma}=\Gamma\oplus\Z \delta$ with
$\deg e_v t^r = \deg e_{vr} = \alpha_v + r\delta$,
$\deg f_v t^r = \deg f_{vr} = -\alpha_v + r\delta$,
$\deg h_v t^r = \deg h_{vr} = r\delta$
and $\deg c=0$,
and the set of roots for either algebra is
\[
\hat\Delta = \{ \alpha+r\delta \mid \alpha\in\Delta, r\in\Z \} \cup \{ r\delta \mid 0\neq r\in Z \}.
\]
The real roots are $\alpha+r\delta$ with $\alpha$ real.
If $\gf$ is of finite type, then $\LL\gf$ is the corresponding
affine Lie algebra, and if $\gf$ is of affine type, then
$\LL\gf$ is a toroidal algebra.

The Grothendieck group $K_0(\Coh\XX)$ was computed by Geigle and Lenzing,
and following Schiffmann~\cite{Schiffmann} it can be
identified with $\hat{\Gamma}$, with
\begin{equation}
\label{e:kzeroelts}
[\OO(r\oc)] = \alpha_* + r\delta,
\quad
[S_a] = \delta,
\quad
[S_{ij}] = \begin{cases}
\alpha_{ij} & (j\neq 0) \\
\delta - \sum_{\ell=1}^{w_i-1} \alpha_{i\ell} & (j=0).
\end{cases}
\end{equation}
The \emph{type} of a sheaf is the corresponding element of $\hat \Gamma$.
The symmetric bilinear form $(-,-)$ on $\Gamma$ extends to $\hat \Gamma$ by defining $(\delta,-)=0$,
and it corresponds to the symmetrization of
the Euler form
\[
\langle [X], [Y] \rangle = \dim\Hom(X,Y) - \dim \Ext^1(X,Y)
\]
on $K_0(\Coh\XX)$.
Now $K_0(\Coh\XX)$ is partially ordered, with the positive cone being the classes
of objects in $\Coh\XX$. By (\ref{e:kzeroelts}) the corresponding partial
ordering on $\hat \Gamma$ has as positive cone $\hat \Gamma_+$ the non-negative linear combinations
of the elements $\alpha_*+r\delta$ ($r\in\Z$), $\delta$, $\alpha_{ij}$ and
$\delta - \sum_{\ell=1}^{w_i-1} \alpha_{i\ell}$.
Clearly every root is positive or negative.

\begin{theorem}
\label{t:kaccoh}
If $\XX$ is a weighted projective line and $\phi\in\hat \Gamma$,
there is an indecomposable sheaf in $\Coh\XX$ of type $\phi$ if
and only if $\phi$ is a positive root.
There is a unique indecomposable for a real root, infinitely many for an imaginary root.
\end{theorem}

This is an analogue of Kac's Theorem \cite{Kac1,Kac2,KR} which describes
the possible dimension vectors of indecomposable representations
of quivers.

We remark that there is a complete classification of indecomposables
if $\gf$ is of finite type \cite{GL}, and also
if $\gf$ is of affine type \cite{LM}. The latter
is essentially equivalent to Ringel's classification \cite{Ringelbook}
of representations of tubular algebras.

Lenzing \cite[\S 4.2]{Lenzing} showed that the
category of torsion-free sheaves on $\XX$ is equivalent to the category
of \emph{(quasi) parabolic bundles} on $\PP^1$ of weight type $(D,\mathbf{w})$, that is,
vector bundles $\pi:E\to \PP^1$ equipped with a flag of subspaces
\[
\pi^{-1}(a_i) \supseteq E_{i1} \supseteq \dots \supseteq E_{i,w_i-1}
\]
for each $i$.
This equivalence is not unique, but it can be chosen so that
if $E$ is a parabolic bundle, then $[E] = \dimv E + (\deg E)\delta$.
Here the \emph{dimension vector} of $E$ is
\[
\dimv E = n_* \alpha_*  + \sum_{i=1}^k \sum_{j=1}^{w_i-1} n_{ij} \alpha_{ij} \in \Gamma,
\]
with $n_* = \rank E$ and $n_{ij} = \dim E_{ij}$.
Observe that the dimension vector is necessarily \emph{strict}, meaning
that $n_* \ge n_{i1} \ge n_{i2} \ge \dots \ge n_{i,w_i-1}\ge 0$.
We can now restate Theorem~\ref{t:kaccoh} as follows.

\begin{corollary}
For each $d\in\Z$ there is an indecomposable parabolic bundle of
dimension vector $\alpha\in \Gamma$ and degree $d$ if and only if $\alpha$
is a strict root for $\gf$. There is a unique indecomposable for a real root,
and infinitely many for an imaginary root.
\end{corollary}

In \cite{CBipb} this result is shown to be related to the
existence of matrices in prescribed conjugacy class closures with
product equal to the identity. Using that, in case the
matrices have generic eigenvalues, we gave a partial proof
over the complex field.

Our proof of Theorem \ref{t:kaccoh} uses Hall algebras.
First we need a lemma, which we have observed with C.~Gei{\ss}.
Given a parabolic bundle $E$, the underlying vector bundle on $\PP^1$
decomposes as a direct sum of line bundles of degrees
$n_1\le \dots \le n_r$.
One might call $n_r - n_1$ the \emph{width} of $E$.

\begin{lemma}
\label{l:bounded}
For any $\phi\in\hat \Gamma$ there is a bound, depending only on $\mathbf{w}$ and $\phi$,
of the width of indecomposable parabolic bundles of type $\phi$.
\end{lemma}

Equivalently, for any $\phi,\psi\in\hat \Gamma$ there is a bound on $\dim \Hom(X,Y)$
(and so also on $\dim\Ext^1(X,Y)$) for $X,Y$ indecomposable of types $\phi,\psi$.

\begin{proof}
The argument is the same as \cite[Theorem 1]{Atiyahe}.
Any torsion-free sheaf $E$ has a splitting by rank-one
torsion-free sheaves $(L_1,\dots,L_r)$, meaning that
there is a chain $0=E_0\subset E_1 \subset \dots \subset E_r = E$
and $L_i = E_i/E_{i-1}$.
The \emph{degree} is defined for weighted projective lines
by \cite[Proposition 2.8]{GL}, and using it one may consider
splittings which are \emph{maximal} in the sense
that $L_1$ has maximal possible degree, and amongst
these $L_2$ has maximal possible degree, etc.

By \cite[Corollary 1.8.1]{GL} and the structure of
the ring $S(\mathbf{w},D)$, it is clear that
there is an integer $h$ with $\Hom(L,L')\neq 0$
for any rank one torsion-free sheaves $L,L'$ with $\deg L' - \deg L > h$.
If $(L_1,L_2)$ is a maximal splitting of $E$, then
there there is an exact sequence
\[
\Hom(L_1(\oc),E) \to \Hom(L_1(\oc),L_2) \to \Ext^1(L_1(\oc),L_1).
\]
The right hand space is zero since $\Ext^1(\OO(\oc),\OO)=0$.
If $\deg L_2 - \deg L_1(\oc) > h$,
then the middle space is nonzero, and so $\Hom(L_1(\oc),E)\neq 0$.
Taking the image of such a map, and enlarging it so that the quotient of
$E$ by this subsheaf is torsion-free, one contradicts the maximality of
the splitting. Thus we must have $\deg L_2 - \deg L_1(\oc) \le h$,
giving a bound of the form $\deg L_2 - \deg L_1 \le h'$, for some $h'$.
As in \cite[Lemma 4]{Atiyahe} this gives bounds $\deg L_i - \deg L_{i-1} \le h'$
for any maximal splitting $(L_1,\dots,L_r)$,
so $\deg L_i \le \deg L_1 + (i-1)h'$.

Now suppose $E$ is indecomposable, and let $(L_1,\dots,L_r)$ be a maximal
splitting. We show by induction that $\deg L_i \ge \deg L_1 - (i-1)h''$
where $h''=\delta(\vec{\omega})$ in the notation of \cite{GL}.
For $1<i\le n$, since $E$ is indecomposable we must have
$\Ext^1(E/E_{i-1},E_{i-1})\neq 0$,
so $\Hom(E_{i-1},(E/E_{i-1})(\vec{\omega}))\neq 0$ by Serre duality,
and hence $\Hom(L_j(-\vec{\omega}),(E/E_{i-1}))\neq 0$ for some $j< i$.
This implies that $E/E_{i-1}$ has a subsheaf of degree at least $\deg L_j-h''$,
so by maximality $\deg L_i \ge \deg L_j-h'' \ge \deg L_1 - (i-1)h''$
by induction.

The assertion follows.
\end{proof}

For an alternative approach see \cite[Theorem 2.9]{LdlP}.

\section{Hall algebras}
Let $K$ be a finite field and let
$\RR$ be a triangulated $K$-category which is \emph{2-periodic},
meaning that the shift functor $T$ satisfies $T^2=1$.
There is a bilinear form on $K_0(\RR)$,
\[
\langle [X],[Y] \rangle = \dim \Hom(X,Y) - \dim \Hom(X,TY),
\]
and let $(-,-)$ be its symmetrization.
Let $\ind \RR$ be a set of representatives of the isomorphism classes
of indecomposable objects in $\RR$.
Assume that $\RR$ is \emph{finitary}, meaning that
it has finite Hom spaces and
$\{ X \in\ind\RR \mid [X] = \phi \}$ is finite for all $\phi\in K_0(\RR)$.
For $X\in\ind\RR$, define
$d(X) = \dim (\End(X)/\rad\End(X))$,
and assume for simplicity that $K_0(\RR)$ is torsion-free, generated
by indecomposables with $d(X)=1$, and that $[X]$ is divisible in
$K_0(X)$ by $d(X)$ for all $X\in\ind\RR$.
Define
\[
F^Z_{XY} = |\{\text{triangles $Y\to Z\to X\to$}\} / \Aut(X)\times\Aut(Y) |.
\]
Let $\Lambda$ be a commutative ring.
Assuming that $|K|=1$ in $\Lambda$, Peng and Xiao \cite{PX,Hubery} proved that
\[
L_\Lambda(\RR) = (\Lambda \otimes_\Z K_0(\RR)) \; \oplus \bigoplus_{X\in \ind \RR} \Lambda u_X
\]
becomes a Lie algebra over $\Lambda$ with bracket
\[
[u_X,u_Y] = \begin{cases}
{\displaystyle \sum_{Z\in\ind\RR} (F^Z_{XY} - F^Z_{YX}) u_Z } & (X\not\cong TY) \\[15pt]
{\displaystyle 1 \otimes \frac{[X]}{d(X)}}                    & (X\cong TY)
\end{cases}
\]
and $[1\otimes \phi,u_X] = -(\phi , [X]) u_X$ and
$[1\otimes \phi,1\otimes \psi] = 0$ for $\phi,\psi\in K_0(\RR)$.

We now consider weighted projective lines over finite fields,
in the case when the marked points are all defined over the finite field.
The category $\Coh\XX$ is still defined and well-behaved,
see see \cite{Lenzing} or \cite{Schiffmann}.
Schiffmann \cite{Schiffmann} has considered its Hall algebra,
and related it to a quantum group for the positive part of $\LL\gf$.
To apply the construction of Peng and Xiao one uses the quotient category
\[
\RR_\XX = D^b(\Coh\XX)/(T^2),
\]
called the \emph{root category}, whose objects are representatives of the orbits
of $T^2$ on $D^b(\Coh\XX)$, and with
\[
\Hom_{\RR_\XX} (X,Y) = \bigoplus_{n\in \Z} \Hom_{D^b(\Coh\XX)} (X, T^{2n} Y).
\]
This is known to be a 2-periodic triangulated category by \cite[Lemma 2.3]{PXroot}. (See
also \cite[\S3]{LP} for the transition from hereditary algebras to hereditary
abelian categories.)
Since $\Coh\XX$ is hereditary, the indecomposable objects in $D^b(\Coh\XX)$
are the shifts of the indecomposables in $\Coh\XX$, and hence
\[
\ind \RR_\XX = (\ind \Coh\XX) \cup \{ TY \mid Y\in \ind \Coh\XX \}.
\]
Recall that any triangle $X\to Y\to Z\to $ can be rotated to give a triangle $Y\to Z\to TX\to$.
Any triangle $X\to Y\to Z\to$ in $\RR_\XX$ with $X,Y,Z$ indecomposable can be
rotated sufficiently so that $X$ and $Z$ are in $\Coh\XX$, and in this case $Y$ must also be,
and then such triangles are in 1-1 correspondence with short exact sequences $0\to X\to Y\to Z\to 0$.

Assuming that the base field $K$ is finite and $|K|=1$ in $\Lambda$,
the construction of Peng and Xiao
gives a Lie algebra $L_\Lambda(\RR_\XX)$
with triangular decomposition
\[
L_\Lambda(\RR_\XX) =
\biggl(\bigoplus_{X\in \ind \Coh\XX} \Lambda u_X \biggr)
\oplus
(\Lambda \otimes_\Z \hat \Gamma)
\oplus
\biggl(\bigoplus_{Y\in \ind \Coh\XX} \Lambda u_{TY} \biggr).
\]
We define $b_X$ for $X\in\ind\RR_\XX$ by
$b_Y = u_Y$ and $b_{TY} = -u_{TY}$
for $Y\in\ind\Coh\XX$.
If $S$ is a simple sheaf, we extend the notation $S[r]$
to $r<0$ by defining $S[r] = TY$,
where $Y$ is the unique uniserial sheaf of length~$-r$ with
$\Ext^1(Y,S)\neq 0$, so that $\Hom(S[r],S)\neq 0$.
Let $H_r$ be the set of $X\in \ind \RR_\XX$
of type $r\delta$ and with $\Hom(X,S_{ij})=0$ for all $1\le i\le k$,
$1\le j\le w_i-1$, and set $\mathbf{h}_r = \sum_{X\in H_r} d(X) b_X$.

\begin{theorem}
\label{t:loppelts}
The following elements of $L_\Lambda(\RR_\XX)$ satisfy the
relations \textup{(\ref{e:loopalg})} for $\LL\gf$.
\begin{align*}
e_{vr} &=
\begin{cases}
b_{S_{ij}[rw_i+1]}  & (v = ij) \\
b_{\OO(r\oc)}          & (v = *),
\end{cases}
\quad
f_{vr} =
\begin{cases}
b_{S_{i,j-1}[rw_i-1]} & (v = ij) \\
b_{T\OO(-r\oc)}          & (v = *),
\end{cases}
\\
c &= - 1\otimes\delta,
\quad
h_{vr} =
\begin{cases}
- 1\otimes \alpha_v                           & (r=0) \\
b_{S_{ij}[rw_i]} - b_{S_{i,j-1}[rw_i]}  & (r\neq 0,\ v=ij) \\
\mathbf{h}_r                                      & (r\neq 0,\ v=*).
\end{cases}
\end{align*}
\end{theorem}

See also \cite{LP}, where elliptic Lie algebra generators are found in
$L_\Lambda(\RR_\XX)$ for $\gf$ of affine type,
\cite{Schiffmann}, where the Hall algebra of $\Coh\XX$ is
considered, and \cite{Kapranov},
where doubled Hall algebras are considered.

\section{Proof of Theorem \ref{t:loppelts}}

\begin{lemma}
\label{l:coeffone}
If $0\to X\to Y\to Z\to 0$ is a short exact sequence of indecomposable
finite-length sheaves, then up to automorphisms of any two of $X,Y,Z$, any other
exact sequence with the same terms is equivalent to this one.
\end{lemma}

\begin{proof}
Since $Y$ is uniserial, it has a unique subsheaf $Y'$
isomorphic to $X$, from which it is clear that
there is a unique sequence up to the
action of $\Aut(X)\times\Aut(Z)$.
For the action of $\Aut(X)\times\Aut(Y)$, say,
we reduce to the case where $X,Y,Z$ are
finite-dimensional modules
for a finite-dimensional serial algebra,
and we may assume that $Y$ is projective.
Then any two epimorphisms $Y\to Z$ are
equivalent via an element of $\Aut(Y)$,
and the result follows.
\end{proof}

\begin{lemma}
\label{l:tsij}
$TS_{ij}[r] = S_{i,j-r}[-r]$
where the subscript $j-r$ is computed modulo $w_i$.
\end{lemma}

\begin{proof}
Clear.
\end{proof}

\begin{lemma}
\label{l:bsbs}
One has
\[
[b_{S_{ij}[r]},b_{S_{ik}[s]}] =
\begin{cases}
\kdelta_{j-r,k}b_{S_{ij}[r+s]} - \kdelta_{j,k-s} b_{S_{ik}[r+s]} & (r+s\neq 0) \\
-\kdelta_{j-r,k} \otimes [S_{ij}[r]] & (r+s=0),
\end{cases}
\]
where the subscripts $j-r$ and $k-s$ are computed modulo $w_i$.
\end{lemma}

\begin{proof}
If $r,s>0$, then one gets a positive contribution of $u_X$
for short exact sequences $0\to S_{ik}[s] \to X \to S_{ij}[r]\to 0$,
and a negative contribution for short exact sequence
$0\to S_{ij}[r] \to X \to S_{ik}[s]\to 0$. The condition for
the existence of nonsplit sequences is given by the $\kdelta$'s.
In each case there is a unique possible middle term, and
the coefficient is 1 by Lemma \ref{l:coeffone}.

If $r,s<0$ the argument is similar.

If $r>0,s<0$, one gets a contribution of $u_X$
for $X$ in a triangle
$S_{ik}[s]\to X\to S_{ij}[r]\to$ or $S_{ij}[r]\to X\to S_{ik}[s]\to$.
Rotating, these become triangles $X\to S_{ij}[r]\to S_{i,k-s}[-s]\to$
and $S_{i,k-s}[-s]\to S_{ij}[r]\to X\to$.
Suppose that $r\ge -s$ (the reverse is similar).
Then $X$ must be a sheaf in both cases,
corresponding to short exact sequences
$0\to X\to S_{ij}[r]\to S_{i,k-s}[-s]\to 0$
and $0\to S_{i,k-s}[-s]\to S_{ij}[r]\to X\to 0$.
The existence of such sequences is given by the $\kdelta$'s,
and in each case there is a unique possible $X$.
\end{proof}

\begin{lemma}
\label{l:osijsext}
There is a short exact sequence
$0\to \OO(r\oc)\to X\to S_{ij}[s]\to 0$
with $X$ indecomposable if and only if $j \equiv s \, (\modu \, w_i)$,
and then $X\cong \OO(r\oc+s\ox_i)$.
\end{lemma}

\begin{proof}
If $X$ is indecomposable it is of the form $\OO(\ox)$ for some $\ox$,
and by considering the type, one must have $\ox=r\oc+s\ox_i$. Now
since there is a nonzero homomorphism $\OO(\ox)\to S_{ij}$, one has
$j \equiv s \, (\modu \, w_i)$.
\end{proof}

\begin{lemma}
\label{l:bdbd}
If $X,Y\in\ind\RR_\XX$ and $[X]=r\delta$, $[Y]=s\delta$ then
$[b_X,b_Y]=0$ if $X\not\cong TY$.
\end{lemma}

\begin{proof}
To have any chance of $[b_X,b_Y]$ being nonzero,
the simple sheaves involved in $X$ and $Y$ must all
be of the form $S_a$ or must all be of the form $S_{ij}$
for fixed $i$. The latter case follows from Lemma~\ref{l:bsbs}.
The former case is analogous.
\end{proof}

\begin{lemma}
\label{l:th}
$H_{-r} = \{ TY \mid Y\in H_r \}$.
\end{lemma}

\begin{proof}
Clear.
\end{proof}

\begin{lemma}
\label{l:sumhr}
$\sum_{X\in H_r} d(X) = 2$ in $\Lambda$.
\end{lemma}

\begin{proof}
We may assume that $r>0$. The restriction
$\Hom(X,S_{ij})=0$ for all $1\le i\le k$,
$1\le j\le w_i-1$, ensures that the marked points
can each contribute at most one indecomposable.
Thus this is a question about torsion sheaves on $\PP^1$.
The point at infinity contributes one
indecomposable sheaf, and the rest correspond to
indecomposable $r$-dimensional modules for
the polynomial ring $K[x]$.
Now absolutely indecomposable modules are given by
Jordan blocks, so the number is equal to the size of the field,
and as this is equal to 1 in $\Lambda$,
formula ($\alpha$) on page 91 of \cite{Kac2}
gives the result.
\end{proof}

We now verify that the elements of Theorem~\ref{t:loppelts}
satisfy the relations \textup{(\ref{e:loopalg})} for $\LL\gf$.
The arguments are all standard in the theory of Hall algebras.

\medskip
\noindent
(i)
$c$ central. This is clear since $(\delta,-)=0$.

\medskip
\noindent
(ii)
$[e_{vr},e_{vs}] = 0$.
\begin{itemize}
\item[(a)]%
If $v=ij$ this follows from Lemma~\ref{l:bsbs}.

\item[(b)]%
If $v=*$ we want $[u_{\OO(r\oc)},u_{\OO(s\oc)}]=0$.
The sheaves $\OO(r\oc)$ all lie in a subcategory of $\Coh\XX$ which
is equivalent to $\Coh\PP^1$. In any extension, the middle term lives
in this category $\Coh\PP^1$, but here the indecomposables are all line bundles.
\end{itemize}

\medskip
\noindent
(iii)
$[f_{vr},f_{vs}] = 0$. Similar to (ii).

\medskip
\noindent
(iv)
$[h_{ur},h_{vs}]=ra_{uv}\,\kdelta_{r+s,0}\,c$.
Expanding the left hand side, observe that every $u_X$ which
occurs has $[X]=r\delta$ or $s\delta$,
so in the radical of the symmetric bilinear form. Thus by Lemma~\ref{l:bdbd},
the only way to
not get zero is if $h_{ur}$ involves a $u_X$ and $h_{vs}$ involves the corresponding $u_{TX}$.
Thus the only possibilities are $[h_{ur},h_{v,-r}]$ with $r\neq 0$.
By symmetry we may assume that $r>0$.
\begin{itemize}
\item[(a)]%
By Lemmas~\ref{l:th} and \ref{l:sumhr} we have
\begin{align*}
[h_{*,r}, h_{*,-r}]
&= \sum_{X,Y\in H_r} d(X) d(Y) [b_X,b_{TY}]
\\
&= - \sum_{X\in H_r} d(X)^2 [u_X,u_{TX}]
\\
&= - \sum_{X\in H_r} d(X)^2 1\otimes [X]/d(X)
\\
&= - 1\otimes \sum_{X\in H_r} d(X) [X]
\\
&= -1\otimes r\delta \sum_{X\in H_r} d(X)
\\
&= 2r(-1\otimes \delta)
= 2rc.
\end{align*}

\item[(b)]%
$[h_{ij,r},h_{ij,-r}] =
[b_{S_{ij}[rw_i]}-b_{S_{i,j-1}[rw_i]},b_{S_{ij}[-rw_i]}-b_{S_{i,j-1}[-rw_i]}].
$
Expanding this, the cross terms vanish by the argument above, giving
\[
[b_{S_{ij}[rw_i]},b_{S_{ij}[-rw_i]}]
+
[b_{S_{i,j-1}[rw_i]},b_{S_{i,j-1}[-rw_i]}]
= -2\otimes r\delta = 2rc.
\]

\item[(c)]%
$[h_{ij,r},h_{k\ell,-r}]$ can only be nonzero, by the argument above, if $k=i$
and $\ell = j$ or $j\pm 1$. If $\ell=j\pm 1$, then one gets a cross term,
so the result is $-rc$.

\item[(d)]%
For $[h_{*,r},h_{ij,-r}]$, the only nonzero term which might
occur comes from $S_{i0}[rw_i]\in H_r$, giving
$[b_{S_{i0}[rw_i]},-b_{S_{i0}[-rw_i]}]$ provided that $j=1$.
This gives~$-rc$.
\end{itemize}

\medskip
\noindent
(v)
$[e_{ur}, f_{vs}] = \kdelta_{uv}\, \left(h_{v,r+s} + r\,\kdelta_{r+s,0}\,c\right)$.

\begin{itemize}
\item[(a)]%
For $[e_{ij,r},f_{k\ell,s}]$, if $r+s=0$ then
\begin{align*}
[e_{ij,r},f_{k\ell,s}]
&= [b_{S_{ij}[rw_i+1]},b_{S_{k,\ell-1}[sw_k-1]}]
\\
&= - \kdelta_{ik} \kdelta_{j-(rw_i+1),\ell-1} \otimes [S_{ij}[rw_i+1]]
\\
&= - \kdelta_{ik} \kdelta_{j\ell} \otimes (\alpha_{ij}+r\delta)
=\kdelta_{ik} \kdelta_{j\ell} (h_{ij,0} + rc),
\end{align*}
and if $r+s\neq 0$ then
\begin{align*}
[e_{ij,r},f_{k\ell,s}]
&= [b_{S_{ij}[rw_i+1]},b_{S_{k,\ell-1}[sw_k-1]}]
\\
&=\kdelta_{ik} \kdelta_{j\ell} \left( b_{S_{ij}[(r+s)w_i]} - b_{S_{i,j-1}[(r+s)w_i]} \right)
\\
&= \kdelta_{ik} \kdelta_{j\ell} h_{ij,r+s}.
\end{align*}

\item[(b)]%
For $[e_{*,r},f_{*,s}]$, if $r+s=0$ then
\begin{align*}
[e_{*,r},f_{*,s}]
&= - [u_{\OO(r\oc)},u_{T\OO(-s\oc)}] =
- 1\otimes [\OO(r\oc)]
\\
&= - 1\otimes (\alpha_* + r\delta) = h_{*,0} + rc,
\end{align*}
so suppose that $r+s\neq 0$. In computing
$[e_{*,r},f_{*,s}] = - [u_{\OO(r\oc)},u_{T\OO(-s\oc)}]$,
one gets a negative contribution of $u_X$ for triangles
$T\OO(-s\oc)\to X\to \OO(r\oc) \to$,
which is only possible when $X = TY$ with $Y$ a uniserial sheaf,
and a positive contribution for triangles
$\OO(r\oc)\to X\to T\OO(-s\oc)\to$, which is possible
for $X = Y$, a uniserial sheaf.
Thus one gets a positive contribution of $b_X$ in each case.
In computing the coefficients, one may apply a shift to the triangles,
so one sees that the answer only depends on $r,s$ through their sum $t=r+s$.
Thus one gets contributions for exact sequences
$0\to \OO(t\oc) \to \OO \to Y \to 0$
and
$0 \to \OO(-t\oc) \to \OO \to Y \to 0$.
Assuming that $t>0$ (the case $t<0$ is similar),
only the latter are involved.
The possible $Y$ are those in $H_t$, and for such $Y$,
if $S$ is the simple in its top, and $d=d(Y)=d(S)$,
then there are $t/d$ copies of $S$ involved in $Y$.
Now $\Hom(\OO,Y)$ has dimension $t$, and the
non-epimorphisms give a subspace of dimension $t-d$.
Thus the number of exact sequences is
\[
(q-1)(q^t - q^{t-d}).
\]
Factoring out by the automorphisms of $\OO(-t\oc)$ and $\OO$,
which act freely, one gets
\[
\frac{q^t - q^{t-d}}{q-1} = q^{t-d} \, \frac{q^d-1}{q-1}.
\]
In $\Lambda$ this is $d$,
so $\sum_{Y\in H_t} d(Y) b_Y = \mathbf{h}_t = h_{*,t}$.

\item[(c)]%
For $[e_{*,r},f_{ij,s}]$,
one gets contributions from triangles
$S_{i,j-1}[sw_i-1] \to X \to \OO(r\oc) \to$
and $\OO(r\oc) \to X \to S_{i,j-1}[sw_i-1] \to$.
Rotating, the first becomes
$X \to \OO(r\oc) \to S_{ij}[-sw_i+1] \to$
by Lemma \ref{l:tsij}.
Now there can be nonzero homomorphisms from $\OO(r\oc)$ to
$S_{ij}[-sw_i+1]$ only if the latter is a sheaf, but then
there are no epimorphisms since $j\neq 0$.
The second becomes
$X \to S_{i,j-1}[sw_i-1] \to T\OO(r\oc) \to $
and there can only be nonzero homomorphisms from $S_{i,j-1}[sw_i-1]$ to $T\OO(r\oc)$
if $S_{i,j-1}[sw_i-1]$ is a sheaf. Thus
one deals with short exact sequences
$0\to \OO(r\oc) \to X \to S_{i,j-1}[sw_i-1] \to 0$.
Since $X$ is indecomposable, it must be a torsion-free sheaf.
Now if $f$ is the morphism $X \to S_{i,j-1}[sw_i-1]$ and $S$ is
the socle of $S_{i,j-1}[sw_i-1]$, then $f^{-1}(S)$ must also
be torsion-free. But the sequence $0\to\OO(r\oc)\to f^{-1}(S)\to S\to 0$
splits since $S\cong S_{i,j+1}$.

\item[(d)]%
$[e_{ij,r},f_{*,s}]$ is similar to (c).
\end{itemize}

\medskip
\noindent
(vi)
$[h_{ur},e_{vs}] = a_{uv} e_{v,r+s}$.
If $r=0$ then
\[
[h_{ur},e_{vs}] = [-1\otimes\alpha_u,e_{vs}] = (\alpha_u, \alpha_v+s\delta)e_{vs}
\]
as required,
so suppose $r\neq 0$. We assume that $r>0$. (The case $r<0$ is similar.)

\begin{itemize}
\item[(a)]%
$[h_{ij,r},e_{k\ell,s}] = [b_{S_{ij}[rw_i]}-b_{S_{i,j-1}[rw_i]},b_{S_{k\ell}[sw_k+1]}]$,
and Lemma~\ref{l:bsbs} gives the result.

\item[(b)]%
$[h_{ij,r},e_{*,s}] = [b_{S_{ij}[rw_i]}-b_{S_{i,j-1}[rw_i]},b_{\OO(s\oc)}]$.
In expanding, one gets contributions $u_X$ only for short exact
sequences with middle term $X$ and end terms the sheaves in the
expression. By the argument in (v)(c),
the only possible extension with indecomposable middle term is
$0\to \OO(s\oc)\to X\to S_{i0}[rw_i]\to 0$, and then $X \cong \OO((r+s)\oc)$.
There is only one such extension, modulo automorphisms,
giving $[h_{ij,r},e_{*,s}] = -b_{\OO((r+s)\oc)} = -e_{*,r+s}$.

\item[(c)]%
$[h_{*,r},e_{ij,s}] = \sum_{X\in H_r} d(X) [b_X,b_{S_{ij}[sw_i+1]}]$.
One gets a contribution of $u_Y$ for triangles
$S_{ij}[sw_i+1]\to Y\to X\to$ and
$X\to Y\to S_{ij}[sw_i+1]\to$.

If $s\ge 0$ these correspond to short exact sequences
$0\to S_{ij}[sw_i+1]\to Y\to X\to 0$ and
$0\to X\to Y\to S_{ij}[sw_i+1]\to 0$.
For the first, there are
no indecomposable $Y$, and for the second
there is only an exact sequence with $Y$ indecomposable if $j=1$ and $X\cong S_{i0}[rw_i]$,
and then $[h_{*,r},e_{ij,s}] = -u_{S_{ij}[rw_i+sw_i+1]} = -e_{ij,r+s}$.

If $s<0$ and $r+s \ge 0$, the triangles correspond to short exact sequences
$0\to Y\to X\to S_{i,j-1}[-sw_i-1]\to 0$ and
$0\to S_{i,j-1}[-sw_i-1]\to X\to Y\to 0$.
and the only possibility is $j=1$ and $Y \cong S_{i,w_i-1}[rw_i+sw_i+1]$
in the first of these, so again
$[h_{*,r},e_{ij,s}] = -u_{S_{ij}[rw_i+sw_i+1]} = -e_{ij,r+s}$.

If $r+s<0$, the triangles correspond to short exact sequences
$0\to X\to S_{i,j-1}[-sw_i-1]\to TY\to 0$ and
$0\to TY\to S_{i,j-1}[-sw_i-1]\to X\to 0$, and the only possibility
is $j=1$ and $TY \cong S_{i0}[-rw_i-sw_i-1]$, and again
$[h_{*,r},e_{ij,s}] = -e_{ij,r+s}$.

\item[(d)]%
$[h_{*,r},e_{*,s}] = \sum_{X\in H_r} d(X) [b_X,b_{\OO(s\oc)}]$.
Computing the brackets on the right hand side, one gets a positive
contribution of $u_Y$ for triangles $\OO(s\oc)\to Y\to X\to$,
and a negative contribution for triangles $X\to Y\to \OO(s\oc) \to$.
In the first case $Y$ must be a sheaf. In the second it
must also be a sheaf, but there are no nonsplit extensions.
Consider exact sequences $0\to\OO(s\oc)\to Y\to X\to 0$.
The only possible $Y$ is $\OO((r+s)\oc)$, and the number
of sequences modulo automorphisms of $\OO(s\oc)$ and $X$ is 1.
Thus $\sum_{X\in H_r} d(X) u_{\OO((r+s)\oc)} = 2e_{*,r+s}$
by Lemma~\ref{l:sumhr}.

\end{itemize}

\medskip
\noindent
(vii)
$[h_{ur},f_{vs}] = -a_{uv} f_{v,r+s}$.
Similar to (vi).

\medskip
\noindent
(viii)
$(\ad e_{u0})^{1-a_{uv}}(e_{vs}) = 0$ for $u\neq v$.
\begin{itemize}
\item[(a)]%
$[e_{ij,0},e_{k\ell,s}]=0$ for $k\neq i$ or $\ell\neq j\pm 1$ by Lemma \ref{l:bsbs}.

\item[(b)]%
$[e_{*,0},e_{ij,s}]=0$ for $j>1$.
One gets a contribution
of $u_X$ for sheaves belonging to short exact sequences
$0\to\OO\to X\to S_{ij}[sw_i+1]\to 0$.
Now the epimorphism $S_{ij}[sw_i+1]\to S_{ij}$ induces
an epimorphism $X\to S_{ij}$. If $L$ is its kernel, then $L$ is
an extension of $\OO$ by $S_{i,j-1}[sw_i]$, so $L\cong \OO(s\oc)$.
But there is no nonsplit extension $0\to \OO(s\oc)\to X\to S_{ij}\to 0$
for $j>1$, so $X$ must decompose.

\item[(c)]%
$[e_{ij,0},[e_{ij,0},e_{i\ell,s}]]=0$ for $\ell=j\pm 1$ by Lemma \ref{l:bsbs}.

\item[(d)]%
$[e_{i1,0},[e_{i1,0},e_{*,s}]]=0$.
Computing $[e_{i1,0},e_{*,s}]$, one gets a contribution of
$u_X$ for short exact sequences $0\to\OO(s\oc)\to X\to S_{i1}\to 0$,
and the only possibility is $X\cong \OO(s\oc+\ox_i)$.
Then, computing $[e_{i1,0},[e_{i1,0},e_{*,s}]]$, one
gets a contribution of $u_Y$ for short exact sequences
$0\to \OO(s\oc+\ox_i)\to Y\to S_{ij}\to 0$, but there
are no nonsplit extensions.

\item[(e)]%
$[e_{*,0},[e_{*,0},e_{i1,s}]]=0$.
Computing $[e_{*,0},e_{i1,s}]$, one gets a contribution
of $u_X$ for short exact sequences
$0\to \OO\to X\to S_{i1}[sw_i+1]\to 0$,
and then one gets a contribution to
$[e_{*,0},[e_{*,0},e_{i1,s}]]$ of $u_Y$ for
short exact sequences
$0\to X\to Y\to S_{i1}[sw_i+1]\to 0$.
Now by the theory of perpendicular categories \cite{GLperp},
all of these sheaves belong to a subcategory of
$\Coh\XX$ corresponding to coherent sheaves on a
weighted projective line with only one marked point, $a_i$,
and for this subcategory it is known by the work of Geigle and Lenzing
\cite{GL} that all indecomposable torsion-free sheaves have rank 1.
But $Y$ would have to have rank 2.

\end{itemize}

\medskip
\noindent
(ix)
$(\ad f_{u0})^{1-a_{uv}}(f_{vs}) = 0$ for $u\neq v$.
Similar to (viii).

\section{Proof of Theorem \ref{t:kaccoh}}
Let $\Phi$ be an additive group,
$(-,-):\Phi\times \Phi\to \Z$ a symmetric bilinear form,
and let $\alpha\in \Phi$ satisfy $(\alpha,\alpha)=2$.
One of the standard arguments in Lie theory
shows that if $L$ is a $\Phi$-graded complex Lie algebra,
$e\in L_\alpha$, $f\in L_{-\alpha}$ and $h=[e,f]$ have the property
that $\ad e$ and $\ad f$ are locally nilpotent and
$\ad h$ acts on any $L_\psi$ as multiplication by $(\alpha,\psi)$,
then $\dim L_\phi = \dim L_{\phi - (\alpha,\phi)\alpha}$
for any $\phi\in \Phi$.
Namely, the operator
$\theta = \exp(\ad e) \exp(-\ad f) \exp(\ad e)$
is defined, and $\theta(h)=-h$.
If $x\in L_\phi$, we can write $\theta(x) = \sum_{r\in \Z} y_r$ with $y_r \in L_{\phi+r\alpha}$,
and
\[
\sum_{r\in\Z} (\alpha,\phi) y_r = \theta([h,x]) = [\theta(h),\theta(x)]
= [-h,\theta(x)] = \sum_{r\in\Z} - (\alpha,\phi+r\alpha) y_r.
\]
Thus, for all $r$ either $y_r=0$ or $(\alpha,\phi) = -(\alpha,\phi+r\alpha)$,
so $r = -(\alpha,\phi)$.
Thus, if $x\neq 0$, $(\alpha,\phi)$ must be an integer,
and $\theta(x)\in L_{\phi-(\alpha,\phi)\alpha}$.
Thus $\theta(L_\phi)\subseteq L_{\phi-(\alpha,\phi)\alpha}$.
Similarly $\theta^{-1}(L_{\phi-(\alpha,\phi)\alpha}) \subseteq L_\phi$.
This argument uses in several places that the base field has
characteristic zero, but clearly it gives the following.

\begin{lemma}
\label{l:ifnilpeq}
Given a function $\nu:\Phi\to \N$ and $\phi\in \Phi$, there is some $\ell_0>0$
with the following property.
If $L$ is a $\Phi$-graded Lie algebra over a field of characteristic $\ell\ge \ell_0$,
and $e\in L_\alpha$, $f\in L_{-\alpha}$ and $h=[e,f]$ have the property that
\[
(\ad e)^{\nu(\psi)}(x) = 0,
\quad
(\ad f)^{\nu(\psi)}(x) = 0,
\quad
(\ad h)(x) = (\alpha,\psi)x
\]
for all $\psi\in \Phi$ and $x\in L_\psi$,
then $\dim L_\phi = \dim L_{\phi - (\alpha,\phi)\alpha}$
\end{lemma}

We are going to apply this lemma to Lie algebras of the form $L=L_\Lambda(\RR_\XX)$.
They are graded by $\Phi=\hat \Gamma$ with $L_0 = \Lambda\otimes_\Z \hat \Gamma$,
and $u_X\in L_\psi$ and $u_{TX}\in L_{-\psi}$ for $X$ an indecomposable sheaf of type $\psi$.

Observe that if $\psi\neq 0$, then $\dim L_\psi$ is the number of indecomposable sheaves
of type $\psi$ if $\psi\in\hat\Gamma_+$, is the number of indecomposable sheaves of
type $-\psi$ if $-\psi\in\hat\Gamma_+$, and otherwise zero.

The next lemma will ensure that the nilpotence conditions can be arranged.
Since we have an inequality on the characteristic of the base field for $L$ in
Lemma~\ref{l:ifnilpeq}, it is essential in the next lemma to have a uniform $\nu$,
independent of the finite field $K$.

\begin{lemma}
\label{l:getnu}
Given a weight sequence $\mathbf{w}$ and vertex $v$,
there is a function $\nu:\hat \Gamma\to \N$ such that
for any weighted projective line $\XX$ of type $\mathbf{w}$
over a finite field $K$, the Lie algebra $L=L_\Lambda(\RR_\XX)$ satisfies
\[
(\ad e_{v0})^{\nu(\psi)}(x) = (\ad f_{v0})^{\nu(\psi)}(x) = 0
\]
for all $\psi\in\hat \Gamma$ and $x\in L_\psi$.
\end{lemma}

\begin{proof}
If $X,Y\in\ind \Coh\XX$, $\Ext^1(X,X)=0$,
and $u_Z$ is involved in $(\ad u_X)(u_Y)$, then
$Z$ is the middle term of a nonsplit exact sequence whose
end terms are $X$ and $Y$, so
$\dim \Ext^1(X,Z)+\dim\Ext^1(Z,X)$ is strictly
less than $\dim \Ext^1(X,Y)+\dim\Ext^1(Y,X)$.
Thus $(\ad u_X)^n (u_Y) = 0$ for $n > \dim \Ext^1(X,Y)+\dim\Ext^1(Y,X)$.
The result now follows from Lemma~\ref{l:bounded},
which still holds for $K$ finite, either by inspecting
the argument, or by using the fact that an indecomposable sheaf of type $\phi$
splits over the algebraic closure of $K$ into summands which all have type
$\phi/d$ for a positive integer $d$ dividing $\phi$.
\end{proof}

\begin{lemma}
Suppose given a weight sequence $\mathbf{w}$, vertex $v$, and $0\neq\phi\in \hat\Gamma_+$.
For any prime $p$ there is a power $p^n$
such that if $\XX$ is a weighted projective line of type $\mathbf{w}$
over a finite field $K$ which contains the field with $p^n$ elements,
then the number of indecomposables sheaves of type $\phi$
is the same as the number of type $\pm s_v(\phi)$.
\end{lemma}

\begin{proof}
As explained above, let $\Phi=\hat \Gamma$ and let $\alpha=\alpha_v$.
Let $\nu$ be given by the previous lemma, and $\ell_0$ by Lemma~\ref{l:ifnilpeq}.

Given $p$, choose $n$ so that $p^n-1$ is
divisible by a prime $\ell\ge \ell_0$,
and let $\Lambda$ be a field of characteristic $\ell$.
As $K$ is a finite field containing the field with $p^n$ elements,
we have $|K|=1$ in $\Lambda$, so the Lie algebra $L=L_\Lambda(\RR_\XX)$ is defined.

By Theorem~\ref{t:loppelts}, the elements $e=e_{v0}\in L_\alpha$ and $f=f_{v0}\in L_{-\alpha}$ satisfy
$[e,f]=h_{v0}$, and then by the definition of $L$ we
have $(\ad h_{v0})(x) = (\alpha,\psi)x$ for $x\in L_\psi$.
Thus we are in the setup of Lemma \ref{l:ifnilpeq}, so $\dim L_\phi = \dim L_{s_v(\phi)}$.
This gives the result.
\end{proof}

Now we work over an algebraically closed field.
Given a weighted projective line $\XX$ of weight type $\mathbf{w}$,
and given $\phi\in\hat\Gamma_+$,
one can find an algebraic variety $V$ with the action of a connected
algebraic group $G$, in such a way that the orbits of $G$ on $V$
correspond 1-1 to isomorphism classes of certain sheaves on $\XX$
of type $\phi$, including all indecomposable ones. Moreover,
although the indecomposable sheaves need only form a constructible
subset $V^{\mathrm{ind}}$ of $V$, the subsets $V_d$ of $V$ consisting
of the points whose $G$-orbit has dimension $d$
are locally closed in $V$, and $V_d^{\mathrm{ind}} = V^{\mathrm{ind}}\cap V_d$ is closed in $V_d$,
hence an algebraic variety. One can then
define the \emph{number of parameters} of indecomposable sheaves of type $\phi$ to be
\[
n(\phi) = \max_d \left\{ \dim V_d^{\mathrm{ind}} - d \right\}
\]
and the \emph{number of top-dimensional families} $t(\phi)$ can be defined as
the sum over $d$ of the number of irreducible components
of $V_d^{\mathrm{ind}}$ of dimension $n(\phi)+d$.

For representations of quivers this setup is explained in \cite{KR}. To convert parabolic bundles
to quivers with relations, we used ``squids'' in \cite[Lemma 5.5]{CBipb}.
The essential extra ingredient here is Lemma~\ref{l:bounded}, which ensures that
there is some integer $N$ (depending only on $\mathbf{w}$ and $\phi$), such
any indecomposable parabolic bundle of type $\phi$ satisfies the following condition: (*) if $E$ is
the underlying vector bundle then the dual of the twist of $E$ by $N$ is generated
by global sections. Twisting by $N$, the parabolic bundles satisfying (*) correspond
to representations of the squid satisfying certain conditions (**) mentioned in
\cite[Lemma 5.5]{CBipb}.
Let $R(\alpha)$ denote the space of representations of the corresponding quiver of the
appropriate dimension vector $\alpha$. The representations satisfying the
relations for the squid define a closed subset of $R(\alpha)$, and then the
representations satisfying (**) define an open subset of that, and this is our variety $V$.
An alternative approach would be to use canonical algebras, as in \cite[\S4]{GL}
and \cite[\S3]{LdlP}, and this would allow to handle also the torsion sheaves.
Another approach would be to use ``quot-schemes''.
We omit the details as we only actually need the results for bundles, as the behaviour of the
indecomposable torsion sheaves is trivial.

\begin{lemma}
Suppose given a weight sequence $\mathbf{w}$, vertex $v$, and $0\neq\phi\in \hat\Gamma_+$.
If $\XX$ is a weighted projective line of weight type $\mathbf{w}$
over an algebraically closed field, then
the number of parameters of indecomposable sheaves, and the number of
top-dimensional families for type $\phi$ is
the same as for type $\pm s_v(\phi)$.
In particular, the number of isomorphism classes
of indecomposable sheaves of type $\phi$ (a finite number or $\infty$) is the
same as the number of type $\pm s_v(\phi)$.
\end{lemma}

\begin{proof}
This is essentially the same as in Kac's Theorem \cite{Kac1,Kac2}, see \cite[\S5]{KR}.
We may assume that the point at infinity isn't a marked point.
Letting
\[
T=\Z[x_1,\dots,x_k,\prod_{i<j}(x_i-x_j)^{-1}],
\]
any ring homomorphism
$\theta:T\to K$ defines a weighted projective line $\XX$ of weight type $\mathbf{w}$
over $K$ with marked points $\theta(x_i)$. Moreover the varieties constructed
above are the $K$-points of suitable schemes of finite type over $T$.

Now $\Ker\theta$ is a prime ideal in $T$,
and by constructibility results for dimensions of schemes,
one obtains a maximal ideal $\mathfrak{m}$ lying over $\Ker\theta$
such that the weighted projective lines over $K$ and
over an algebraic closure of the finite field $T/\mathfrak{m}$
have the same numbers of parameters and top-dimensional families for types $\phi$ and $\pm s_v(\phi)$.

This reduces one to the case when $K$ is the algebraic closure
of a finite field. Now if $K_0$ is a finite subfield
containing the marked points, it suffices to show that
over any finite field $K'$ containing $K_0$, the numbers of $K'$-points
of the schemes for types $\phi$ and $s_v(\phi)$ correspond.
This amounts to showing that the numbers
of isomorphism classes of absolutely indecomposable sheaves of types $\phi$ and $s_v(\phi)$
are equal for the corresponding weighted projective line over $K'$.
By an argument involving minimal fields of definition (always containing $K_0$),
it suffices to show that the numbers of indecomposable
sheaves of types $\phi/d$ and $s_v(\phi/d)$ are equal for all $K'$ containing $K_0$
and all positive integers $d$ dividing $\phi$.
This follows from the last
lemma, provided one takes $K_0$ large enough.
\end{proof}

Now let $\phi = \alpha+r\delta\in\hat \Gamma_+$.
If $\alpha=0$ there are infinitely many indecomposables $S_a[r]$
of type $\phi$.
If $\alpha$ is a real root, by a sequence of reflections
one reduces to $\pm\alpha_v + r\delta$, when there
is a unique indecomposable. If $\alpha$ is an imaginary root, one reduces
to $\alpha+r\delta$ with $\alpha$ in the fundamental region,
and there are infinitely many indecomposables by \cite[Lemma 5.6]{CBipb}.
If $\alpha$ is not a root, one reduces to the case
when $\alpha$ is not positive or negative, or has disconnected support,
and there is no indecomposable.
This completes the proof of Theorem~\ref{t:kaccoh}.

\end{document}